\documentclass[a4paper, 10pt, twoside]{article}
\usepackage[english]{babel}
\usepackage{amsmath, amssymb,theorem}
\usepackage[all]{xy}
\usepackage{epsfig,color,rotating}
\xyoption{2cell}
\UseAllTwocells
\UseComputerModernTips
\newcommand{\note}{\textit}
\theoremheaderfont{\normalfont\bfseries}
\theorembodyfont{\normalfont\itshape}
\theoremstyle{change}
\newtheorem{Tm}{\normalfont\scshape Theorem}[section]
\newtheorem{Pp}[Tm]{\normalfont\scshape Proposition}
\newtheorem{Lm}[Tm]{\normalfont\scshape Lemma}
\newtheorem{Cr}[Tm]{\normalfont\scshape Corollary}

\theorembodyfont{\rmfamily}

\newtheorem{Df}[Tm]{\normalfont\scshape Definition}
\newtheorem{Ex}[Tm]{\normalfont\scshape Example}
\newtheorem{Rm}[Tm]{\normalfont\scshape Remark}

\newenvironment{Pf}{{\par\scshape Proof
}}{\hspace*{\fill}{\scshape QED}\\ \vspace*{0.5cm}\par}
\newcommand{\be}{\begin{equation}}
\newcommand{\ee}{\end{equation}}
\makeatletter
\renewcommand{\@makecaption}[2]{%
\vspace{10pt}\hspace{.1\linewidth}
\parbox[l]{.8\linewidth}{\footnotesize{\scshape #1:} #2}%
\par
}
\makeatother

\newcommand{\NN}{\mathbb{N}}

\renewcommand{\phi}{\varphi}

\newcommand{\eps}{\varepsilon}

\newcommand{\cobar}{\smash{\raise.5\baselineskip\hbox{\begin{turn}{180}$B$\end{turn}}}}
\newcommand{\Aut}{\mathrm{Aut}}
\newcommand{\id}{\mathrm{id}}

\newcommand{\clcdot}{\cdot\ldots\cdot}

\title{Families of Hopf algebras of trees \\ and pre-Lie algebras}
\author{Pepijn van der Laan and Ieke Moerdijk}
\begin{document}
\maketitle
\abstract{Using methods from \cite{Moer:Hoa}, we study families of Hopf 
algebra structures on coloured trees.}

\section{Introduction}
Over the past years more and more examples of combinatorial Hopf
algebras appeared in the mathematical literature (cf. Kreimer
\cite{Krei:Hoa},
Connes-Kreimer \cite{ConKr:Hoa}, Loday-Ronco \cite{LodRon:Trees},
Brouder-Frabetti \cite{BrouFrab:Trees}). In most 
cases these Hopf algebras are constructed one at a time. One of the authors
\cite{Moer:Hoa} constructs such Hopf algebras in families stemming
from Hopf operads, rather than as isolated examples. For example, the
Connes-Kreimer Hopf algebra of rooted trees and its planar analogon
(cf. Foissy\cite{Fois:PlantreeII} for a detailed account) are examples
related by a change of operad (see \cite{Moer:Hoa}). The
non-planar version corresponds to a particular coproduct based on the
commutative operad, whereas the planar version corresponds to its
analogue for the associative operad. For the
construction of some other previously studied Hopf algebras of 
trees in the framework of this paper we refer to \cite{Pep:thesis}.  

The paper starts (Section \ref{hoaoninalg}) with the definition of the
initial pair of a commutative algebra together with an $n$-ary map
$(C_n,\lambda_n)$ and applies the arguments of \cite{Moer:Hoa} to
show that there exists a family of Hopf algebra structures on $C_n$. 
Section \ref{inobtree} identifies $C_n$ with the symmetric algebra on rooted
trees with $n$-coloured edges. Consequently the Hopf algebras
constructed in Section \ref{hoaoninalg} are bialgebras of trees. Section
\ref{hoatree} then derives an explicit formula for the corresponding
coproducts in terms of trees. Section \ref{primitive} gives an
explicit description of the Lie algebra 
of primitive elements of the dual Hopf algebra and gives a criterion
when this is the associated Lie algebra of a pre-Lie algebra. Section
\ref{deform} interprets the family of coproducts in terms of deformation
theory. Finally, Section \ref{asscase} sketches the more general framework of
Hopf operads and lists the results obtained when one starts from
associative instead of commutative algebras.

\textbf{Acknowledgements:} During the period in which this article was
written, PvdL was supported by Marie Curie Training Site Fellowship
HPMT-2000-00075 (Centre de Recerca Matem\`atica, Barcelona). 

\section{Hopf algebras structures from initial objects}
\label{hoaoninalg}

For definiteness, we work in the category of vector spaces over a field $k$ of 
characteristic zero. (However, with the exception of Sections
\ref{primitive} an \ref{deform}, our arguments apply in a much  
more general context of modules over an algebra in any symmetric
monoidal additive category.) By algebra we will always mean
associative algebra with unit.

\begin{Df}\label{Ex:treesCnAn}
A \note{commutative $n$-algebra} is a pair 
$(A,\alpha)$, consisting of a  commutative algebra $A$ and a linear map
$\alpha:A^{\otimes n} \longrightarrow A$. A morphism of $n$-algebras
$f:(A,\alpha)\longrightarrow (B,\beta)$ 
is an algebra homomorphism $f:A\longrightarrow B$ 
such that $\beta\circ f^{\otimes n} = f\circ\alpha$. 
We will write $\mathcal{C}_n$ for the category of these $n$-algebras.
For general reasons, this category $\mathcal{C}_n$ has an initial object, 
the free $n$-algebra 
on the empty set of generators. This initial algebra will be denoted by 
$(C_n,\lambda_n)$. It is completely characterised (up to isomorphism) by
the property that any $n$-algebra $(A,\alpha)$ admits a unique
morphism $(C_n,\lambda_n)\longrightarrow (A,\alpha)$. 
We will give an explicit description of this initial algebra in
Proposition \ref{Rm:doehetmetbomen} below.
\end{Df}

This section introduces Hopf algebra structures on $C_n$  of a
specific type. First we need a notation.  
If $(A,\alpha)\in \mathcal{C}_n$ is a commutative $n$-algebra and
$\sigma_1,\sigma_2:A^{\otimes n}\longrightarrow A$ are two linear maps, 
define the linear map
$(\sigma_1,\sigma_2):(A\otimes A)^{\otimes n} \longrightarrow
A\otimes A$ as
\[
(\sigma_1,\sigma_2) = (\sigma_1\otimes \alpha + \alpha \otimes \sigma_2)
\circ \tau.
\]
(Here $\tau:(A^{\otimes 2})^{\otimes n}\longrightarrow
A^{\otimes n}\otimes A^{\otimes n}$ is the isomorphism that
separates the first and second tensor factors from $A^{\otimes 2}$.) 

We are going to use this in the context where $(A,\alpha)$ is the initial $n$-algebra
$(C_n,\lambda_n)$. First note that $(C_n,\lambda_n)$ is an augmented algebra. 
Indeed, the ground field $k$ is naturally an $n$-algebra when
equipped with the zero map $k^{\otimes n}\longrightarrow k$. So by initiality
of $(C_n,\lambda_n)$, there is a unique morphism of $n$-algebras
$\eps:(C_n,\lambda)\longrightarrow (k,0)$.

For any pair of $n$-ary linear maps 
$\sigma_1,\sigma_2:C_n^{\otimes n} \longrightarrow C_n$ as above there
is a unique morphism $\Delta:(C_n,\lambda_n)\longrightarrow (C_n\otimes
C_n,(\sigma_1,\sigma_2))$ in $\mathcal{C}_n$. That
is, a unique algebra morphism $\Delta$ such that the diagram
\[
\xymatrix{ C_n^{\otimes n} \ar[r]^{\lambda}\ar[d]^{\Delta^{\otimes n}} 
& C_n\ar[d]^{\Delta} \\
(C_n\otimes C_n)^{\otimes n}\ar[r]^{(\sigma_1,\sigma_2)}
& C_n\otimes C_n}
\]
commutes. 

\begin{Tm} \label{Tm:genmoerdijk}
Let $n \in \NN$, $\lambda_n$ and $C_n$ be defined as above. Let
$\sigma_i:C_n^{\otimes n}\longrightarrow C_n$ for $i=1,2$ be linear maps. 
If both $\sigma_i$ satisfy
\[
\begin{split}
\epsilon \circ \sigma_i &= \epsilon^{\otimes n},\quad\text{and}\\
\Delta\circ\sigma_i
&= (\sigma_i\otimes\sigma_i)\circ\tau\circ\Delta^{\otimes n}; 
\end{split}
\]
then there exists a unique bialgebra structure on $C_n$
such that $\Delta\circ\lambda = (\sigma_1,\sigma_2)\circ
\Delta^{\otimes n}$ and $\eps\circ\lambda = 0$. 
\end{Tm}

\begin{Pf} 
The proof that this provides a bialgebra structure is completely 
analogous to the case where $n=1$, treated in \cite{Moer:Hoa}. 
\end{Pf}
\begin{Rm}
We will see in the next section that the initial algebra $C_n$
has a natural grading. The coalgebra structure given by a pair of maps
$\sigma_1$ and $\sigma_2$ will respect this grading if these maps $\sigma_1$ and
$\sigma_2$ do. In this case there exists an antipode for the
bialgebra structure. A special case of this occurs in Theorems
\ref{Cr:Hopfspelledout} and \ref{Tm:Coprod} below.
\end{Rm}
\begin{Rm}\label{Rm:simpl}
The algebras $C_n$ together form a simplicial algebra. The
simplicial operations $d_i:C_{n}\longrightarrow C_{n-1}$ and
$s_i:C_{n}\longrightarrow C_{n+1}$ are the algebra homomorphisms
determined by
\[
\begin{split}
d_0(\lambda_n)(x_1,\ldots,x_n) &= \mu(x_1,\lambda_{n-1}(x_2,\ldots,x_n)) \\
d_i(\lambda_n)(x_1,\ldots,x_n) &=
\lambda_n(x_1,\ldots,\mu(x_{i},x_{i+1}),\ldots,x_n)\qquad (i=1,\ldots,n-1)\\
d_n(\lambda_n)(x_1,\ldots,x_n) &= \mu(\lambda(x_1,\ldots,x_{n-1}),x_n)\\
s_i(\lambda_n)(x_1,\ldots,x_n) &=
\lambda_{n+1}(x_1,\ldots,x_{i},1,x_{i+1},\ldots,x_n)\qquad (i=0,\ldots,n),
\end{split}
\]
similar to the formulas for the Hochschild complex.
\end{Rm}
\section{Initial algebras and trees}
\label{inobtree}

\begin{Df}
A \note{rooted tree}  is an isomorphism class $t$ of finite 
partially ordered sets which
\begin{enumerate}
\item have a minimal element $r$ ($\forall x\neq r: r<x$), called
 the \note{root}, and 
\item satisfy the tree condition that 
$(y\neq z) \land (y<x) \land (z<x)$ implies $(y<z)\lor (z<y)$.
\end{enumerate}

In general, we will not be very precise in distinguishing between an isomorphism 
class $t$
and any of the posets which represent it; in particular, we will
often use $t$ to denote to a representing poset, and refer to it as a tree.

The elements of a tree are called \note{vertices}. A pair of vertices
$v<w$ is called an \note{edge} if there is no vertex $x$ such that $v<x<w$.
The number of vertices of a tree $t$ is denoted by $|t|$.
A \note{path} from $x$ to $y$ in a tree is a sequence $(x_i)_i$ of
elements $x=x_n>x_{n-1}>\ldots>x_1>x_0=y$ of maximal length. We will
say that $x$ is above $y$ in a tree if there is a path from $x$ to $y$.
%

A \note{forest} is a finite (possibly empty) multiset (i.e. a set with 
multiplicities) of trees. A subforest of a tree or forest is a subset
of vertices with the induced partial order.

In the sequel we need \note{trees with coloured edges}. These are 
isomorphism classes of posets as above, equipped with a function from 
the set of edges to a fixed set of colours. The isomorphisms are required to 
respect the colours. 
In particular, an \note{$n$-coloured tree} is such a tree whose edges are 
coloured 
by the set $\{1,..,n\}$ of colours. We will write $T_n$ for the vector space 
spanned by 
the set of such $n$-coloured trees. 
\end{Df}

\begin{Pp}\label{Rm:doehetmetbomen}
There is a natural algebra isomorphism between the initial $n$-algebra $C_n$
and the symmetric algebra $S(T_n)$ on the set of $n$-coloured trees.
\end{Pp}

\begin{Pf} 
The symmetric algebra $S(T_n)$ can be identifyed with the vector space 
spanned by the set of $n$-coloured forests, with the unit  represented by the
empty forest and the product by the disjoint union of forests. There is an 
operation
\[
\lambda:S(T_n)^{\otimes n}\longrightarrow S(T_n)
\]
which takes an $n$-tuple of $n$-coloured forests $f_1,\ldots,f_n$, and combines 
them into a single $n$-coloured tree by adding a new root, while connecting this 
new root to each of the roots in the forest $f_i$ by an edge of colour $i$.


This operation makes $S(T_n)$ into an object of the category $\mathcal{C}_n$.  
Since initial 
objects are unique up to isomorphism in any category, it now suffices to prove 
that $S(T_n)$  is  initial in $\mathcal{C}_n$.

To this end, let $(A,\alpha)$ be any object of $\mathcal{C}_n$, where
$\alpha:A^{\otimes n} \longrightarrow A$. Define a morphism 
\[
\phi: (S(T_n),\lambda)\longrightarrow (A,\alpha)
\]
by induction on trees and forests. If $f={t_1\cdot\ldots\cdot t_k}$ is a forest consisting 
of $k$ trees, 
then $\phi(f)=\phi(t_1)\cdot\ldots\cdot \phi(t_k)$, so it suffices to define $\phi$  on 
trees. If $t$
is a tree consisting of a root only, then $\phi(t)=\alpha(1,\ldots,1)$. If $t$ consists 
of a root $r$ onto 
which an $n$-tuple of $n$-coloured forests $f_1,\ldots,f_n$ is
attached by joining the root of each tree in $f_i$ to $r$ via an edge of colour 
$i$, then 
$\phi(t)=\alpha(\phi(f_1),\ldots, \phi(f_n))$. It is straightforward to check that 
$\phi: (S(T_n),\lambda)\longrightarrow (A,\alpha)$ thus defined is indeed a
morphism in $\mathcal{C}_n$, and is the unique such.
\end{Pf}

\section{Hopf algebras of trees}

\label{hoatree}
In this section, we study a particular example of a family of Hopf algebras 
which can be obtained by the general method of Theorem \ref{Tm:genmoerdijk}.

\begin{Ex}\label{Ex:sigma1choice}\label{Ex:deform} 
For any choice of $q_{i\,j}\in k$ for all $j\leq n$, the maps
\[
\sigma_i(t_1,\ldots,t_n)=q_{i\,1}^{|t_1|}\clcdot q_{i\,n}^{|t_n|}
t_1\clcdot t_n \qquad \text{ for } i=1,2
\]
define a bialgebra structure on $S(T_n)$.
Here and in the sequel $|f|$ denotes the number of vertices
in the forest corresponding to $f\in S(T_n)$, while the associative 
multiplication on $S(T_n)$ is denoted by $\cdot$.
\end{Ex}

Below we write $\Delta(f_i) = \sum f_i' \otimes f_i''$ reminiscent of
the form $\Delta$ takes in a basis.

\begin{Tm}\label{Cr:Hopfspelledout}
The symmetric algebra $S(T_n)$ on $n$-coloured trees
has a natural family of graded connected Hopf algebra structures,
indexed by sequences
$(q_{11},\ldots,q_{1n},q_{21},\ldots,q_{2n})\in k^{2n}$. The grading
is with respect to the number of vertices of the trees.
An inductive description of the coproduct is given by
\[
\begin{split}
\Delta(\lambda(f_1,\ldots,f_n)) = &\sum q_{11}^{|f'_1|}\clcdot q_{1n}^{|f'_n|}
\cdot f'_1\clcdot f'_n \otimes \lambda(f''_1,\ldots,f''_n) \\ & + \sum
\lambda(f'_1,\ldots,f'_n)\otimes q_{21}^{|f''_1|}\clcdot q_{2n}^{|f''_n|}
\cdot f''_1\clcdot f''_n,
\end{split}
\]
where $\lambda(f_1,\ldots,f_n)$ is the rooted tree obtained
for $n$ forests $f_1,\ldots,f_n$ by adding a
new root and connecting each of the roots of trees in $f_i$ to the new
root by an edge of colour $i$, and where $|f_i|$ is the number of vertices
in the forest $f_i$.
\end{Tm}
\begin{Pf}
The bialgebra structures are a direct translation of Example
\ref{Ex:deform} to the language of trees of Proposition
\ref{Rm:doehetmetbomen}. The bialgebra $S(T_n)$ is graded connected with respect to
the grading $|.|$. It is well known (cf. Milnor and Moore
\cite{MilnorMoore:Hoa}) that  any graded connected bialgebra admits an
antipode. 
\end{Pf}

We now turn to the question of finding a more direct decription of
these Hopf algebra structures.
Fix $n\in \NN$. For $i=1,2$ and $1\leq j \leq 
n$, let $q_{ij}\in k$, and define for $t_i\in S(T_n)$
\begin{equation}\label{Eq:Convqij}
\sigma_i(t_1,\ldots t_n) = (\prod_j q_{ij}^{|t_j|})\cdot t_1\clcdot t_n.
\end{equation}
Any rooted tree has a natural partial ordering on its vertices in
which the root is the minimal element.  
A \note{subforest} $s$ of a rooted tree $t$ is a subset of the
partially ordered set (representing) $t$ with the induced partial ordering. 
For $v\in s$ we denote by $p_k(v,s,t)$
the number of edges of colour $k$ in the path in $t$ from $v$ to the
root of $t$ that have their  lower vertex in $s^c$. For forests $t$ we
define $p_k(v,s,t)$ as $p_k(v,s\cap t',t')$, where $t'$ is the
connected component of $t$ containing $v$.
There is an easy but useful lemma on the calculus of the $p_k$.

\begin{Lm}
Let $t$ and $s$ be subforests of a forest $u$. Let $v\in s$ and set 
$t' = t\cup v$, $s' = s\cap t'$, $t''= t^c\cup v$ and $s'' = s\cap t''$. Then 
\[
p_k(v,s,u) = p_k(v,s',t') + p_k(v,s'',t''),
\]
where $t',\ t'',\ s'$ and $s''$ are interpreted as subforests of $u$. 
\end{Lm}
\begin{Pf}
The lemma follows at once when we observe that a vertex in the path from 
$v$ to the root in $u$ that is not in $s$ is either in $t'$ or in $t''$.
\end{Pf}
Define for $s\subset t$ a subforest
\begin{equation}
q(s,t):=
\prod_j \left(\prod_{{v\in s}}q_{1\,j}^{p_j(v,s,t)}\cdot
\prod_{v\in s^c}q_{2\,j}^{p_j(v,s^c,t)}\right).
\end{equation}
More intuitively, $q(s,t)$ counts for $v\in s$ the number of edges of colour $j$
in the path from $v$ to the root that have their lower vertex in
$s^c$ and adds a factor $q_{1j}$ for each of these, and $q(s,t)$
counts for $v\in s^c$ the number of edges of colour $j$
in the path from $v$ to the root that have their lower vertex in
$s$ and adds a factor $q_{2,j}$ for each of these.

\begin{Tm}\label{Tm:Coprod}
Let $S(T_n)$ be the symmetric algebra on $n$-coloured trees as in
Theorem \ref{Cr:Hopfspelledout}.  
\begin{enumerate}
\item
For a forest $t\in S(T_n)$ the coproduct defined by
$(q_{11},\ldots,q_{2n})\in k^{2n}$ is given by the formula 
\[
\Delta(t) = \sum_{s\subset t} 
q(s,t) \, s\otimes s^c ,
\]
where the sum is over all subforests $s$ of $t$.
\item
The antipode of the Hopf algebra $S(T_n)$ with the coproduct of part (i) is given by
\[
S(t) = \sum _{k=1}^{|t|}\sum_{\cup_i s_i = t} (-1)^k s_1\clcdot
s_k 
\prod_{1\leq j <k}q(s_j,s_j\cup\ldots\cup s_k),
\]
where we only sum over (ordered) partitions $t = s_1\cup\ldots \cup s_k$ of the
forest $t$ with all forests $s_i$ non-empty.
\end{enumerate}
\end{Tm}
\begin{Pf}
We use induction with respect to the number of applications of
$\lambda$ to show the first result. The formula is trivial for the
empty tree. 
Let $t=\lambda(x_1,\ldots,x_n)$ be a tree and suppose (as
the induction hypothesis) that the formula holds for all 
trees with less then $|t|$ vertices. Since $\Delta$ is an algebra
morphism, it is clear the formula also holds for the forests $x_i$
since these consist of trees with less than $|t|$ vertices. 
Subforests of $t$ are either of the form $ s = \cup_i s_i$, a
(disjoint) union of subforests of the $x_i$, or of the form $s = r
\cup(\cup_i s_i)$, a (disjoint) union of subforests of the $x_i$
together with the root. 
By definition and the induction hypothesis, 
\[
\begin{split}
\Delta(t) &= \sum_{s_i\subset x_i} s_1\clcdot s_n \otimes 
\lambda(s_1^c,\ldots,s_n^c)\cdot
\prod_i q_{1i}^{|s_i|} q(s_i,x_i)\\
&+ \sum_{s_i\subset x_i} \lambda(s_1,\ldots,s_n) \otimes s_1^c \clcdot 
s_n^c\cdot
\prod_i q_{2i}^{|s^c_i|} q(s_i,x_i).
\end{split}
\]
But by the lemma above,
\[
\prod q_{1i}^{|s_i|} q(s_i,x_i) = 
\prod_j \left(\prod_{v\in s}q_{1\,j}^{p_j(v,s,t)}\cdot
\prod_{v\in s^c}q_{2\,j}^{p_j(v,s^c,t)}\right)
\]
for $s=\cup_i s_i = s_1\clcdot s_n$ and $s^c = r\cup(\cup_is_i^c)=
\lambda(s_1^c,\ldots,s_n^c)$; and
\[
\prod q_{2i}^{|s^c_i|} q(s_i,x_i) = 
\prod_j \left(\prod_{v\in s}q_{1\,j}^{p_j(v,s,t)}\cdot
\prod_{v\in s^c}q_{2\,j}^{p_j(v,s^c,t)}\right)
\]
for $s= r\cup(\cup_i s_i) = \lambda(s_1,\ldots,s_n)$ and $s^c =
\cup_is_i=s_1^c\clcdot s_n^c$. Putting these together proves the
formula for the coproduct. 

Let $A$ be a graded connected bialgebra. The augmentation ideal of $A$
is $\bar A = \bigoplus_{n\geq 1} A^n$. 
The antipode on $A$ applied to $x\in \bar A$ is given by
\[
S(x)= \sum_{k=0}^\infty (-1)^{k+1}\mu^{(k)}\circ\bar{\Delta}^{(k)}(x),
\]
where $\bar{\Delta}= \Delta - (\id\otimes 1 + 1 \otimes \id)$, and
$\mu^{(0)}=\id = \bar{\Delta}^{(0)}$, and $\mu^{(k)}:A^{\otimes 
k+1}\longrightarrow A$ and
$\bar{\Delta}^{(k)}:A\longrightarrow A^{\otimes k+1}$ are defined using
(co)associativity for $k>0$ and $\mu^{(0)}=\id = \bar{\Delta}^{(0)}$.
(The sum in the formula for $S(x)$ is of course finite, and stops at $k$ 
for a homogeneous element $x$ of degree $k$.)
As a special case, consider $S(T_n)$ with the coproduct just
described. This proves the result. 
\end{Pf}
\begin{Ex}\label{Ex:Connes-Kreimer}
The Hopf algebra $C_1$ with the
coproduct defined by $q_{11}=1$ and $q_{21} = 0$ is the Connes-Kreimer
Hopf algebra of rooted trees \cite{Krei:Hoa,ConKr:Hoa,Moer:Hoa}.
\end{Ex}

\section{Primitives of the dual}

\label{primitive}
Since $S(T_n)$ is commutative, we know by the Milnor-Moore Theorem
\cite{MilnorMoore:Hoa} that the graded
linear dual $S(T_n)^*$ is the universal enveloping algebra of the Lie
algebra of its primitive elements.  The result below provides an
explicit formula for the Lie bracket on 
the these primitive elements.

\begin{Cr}\label{HnLiefomrula}
Let $S(T_n)$ be a the symmetric algebra on rooted trees with $n$-coloured
edges, and let $\Delta$ be the coproduct defined by
$(q_{11},\ldots,q_{2n})\in k^{2n}$ (cf. Theorem \ref{Tm:Coprod}). The
graded dual $S(T_n)^*$ is the universal enveloping algebra of the Lie
algebra which as a vector space is spanned by elements $D_t$, where
$t$ is a rooted tree in $S(T_n)$. The bracket is given by $[D_s,D_t] =
D_t \bullet D_s-D_s \bullet D_t$, where  
\[
D_t \bullet D_s = \sum_{w} \sum_{s\subset w,\ s^c=t}q(s,w) D_w. 
\]
In this formula, the first sum ranges over all rooted trees $w$, and the
second sum over subtrees of $w$ which are isomorphic to $s$ and whose
complement $s^c$ is isomorphic to $t$.
\end{Cr}

\begin{Pf}
For any cocommutative Hopf algebra we can define an operation $\bullet$ on
the primitive 
elements, such that its commutator is the Lie bracket on primitive
elements: Simply define $\bullet$ as the truncation of the product at
degree $>1$, with respect to the primitive filtration $F$. 
In this case, $F_mC^*_n$ is spanned by the elements
$D_u$ dual to forests $u$ consisting of at most $m$ trees.
The product in $S(T_n)^*$ is determined by the coproduct in
$S(T_n)$. For forest $u$, we can write the multiplication in
$C^*_n$ as
\[
(D_t D_s) (u) = (D_t \otimes D_s) \Delta(u),
\]
thus, for a fixed $u$ we get a contribution $q(u,s)D_u$ for every
subtree isomorphic to $s$ in $u$ with $t$ as complementary forest.
When we then restrict to the primitive part (i.e. the part where $u$
is a tree as well), we conclude that $D_t\bullet D_s$ is given by the
desired formula.
\end{Pf}

Recall (cf. Chapoton and Livernet
\cite{ChapLiv:Prelie}) that a pre-Lie algebra is  vector space $L$
together with a bilinear operation $\bullet$ satisfying the identity
\[
(x\bullet y)\bullet z - x\bullet (y\bullet z) = (x\bullet z)\bullet y
- x\bullet (z\bullet y).
\]
The  the \note{free pre-Lie algebra} 
$L_n$ on $n$ generators is given by the vector space spanned by rooted
trees with vertices labelled by elements of the set $\{1,2,\ldots,n\}$. 
The pre-Lie algebra product is given by grafting trees: For $t$ and
$s$ trees, and $v$ a vertex in $t$, denote by $t\circ_vs$ the tree
obtained from $t$ and $s$ by attaching the root of $s$ to vertex $v$
in $t$ by a new edge. Grafting preserves the labeling of the
vertices. The pre-Lie algebra structure on $L_n$ is given by,
\[
t\bullet s = \sum_{v\in t} t\circ_v s
\]
for trees $s$ and $t$. 

Below we denote by $\chi_S$ the characteristic function of a
subset $S\subset X$ which has value $1$ on $S$ and value $0$ on
$X-S$, and denote the primitive elements of a coalgebra $C$ by $P(C)$.

\begin{Tm}\label{Tm:preLie}
Let $\mathbf{p}\subset\{1,\ldots, n\}$ and define
$q_{1j}=\chi_{\mathbf{p}}(j)$ and $q_{2j}=0$ for $j = 1,\ldots, n$. 
Consider the Hopf algebra structure on the symmetric algebra $S(T_n)$ on rooted
trees with $n$-coloured edges that corresponds to this choice of $q_{ij}$.
\begin{enumerate}
\item The product $D_t\bullet D_s = \sum_{w}\sum_{s\subset w,s^c= t} q(s,w)\, D_w$ of
Corollary \ref{HnLiefomrula} defines a pre-Lie 
algebra structure on the primitive elements $P(C^*_n)$ of $S(T_n)^*$.
\item If $\mathbf{p}=\{1,\ldots,n\}$, then there is a natural inclusion
of this pre-Lie algebra into the free pre-Lie algebra on $n$
generators. The image in the free pre-Lie algebra is spanned by
all sums $\sum_{i\in\mathbf{p}}t_i$, of trees with vertices coloured by
$\mathbf{p}$ that only differ in that the colour of the root of $t_i$
is $i$. 
\end{enumerate}
\end{Tm}
\begin{Pf}
Consider the general formula for $D_t\bullet D_s$ in
Corollary \ref{HnLiefomrula}. Note that for $q_{ij}\in \{0,1\}$ the
coefficients $q(s,w)$ are either 0 or 1. We can be more precise.
Let $w$ and $t$ be trees. A product of subtrees $s=s_1\cdot\ldots\cdot
s_m\subset w$ is $t$-admissible if
$s^c$ contains the root of $w$ while $s$ is grafted onto $t=s^c$ by a
edges of colours $i_1,\ldots,i_m\in \mathbf{p}$ to vertices
$v_1,\ldots,v_m$ respectively each of which is connected to the root
by edges having colours in $\mathbf{p}$. We only use this terminology
for $n=1,2$. Thus, in this case 
we consider $q(s,w)=0$ unless the the corresponding subtree $s$ is $t$-admissible.

The pre-Lie identity follows from
\[
(D_t \bullet D_s)\bullet D_u - D_t \bullet (D_s\bullet D_u) = 
\sum_w\sum_{s\cdot u \subset w} D_w,
\]
where the second sum is over $t$-admissible products $s\cdot u$. This
proves part \textit{(i)} since the expression is symmetric in $s$ and
$u$.

Before we prove \textit{(ii)} we study the pre-Lie algebra structure
$\bullet$ is a bit more detail. Let $m(t,s,w)$ be  the number of
$t$-admissible subtrees $s\subset w$. The operation $\bullet$ is then
given by 
\[
D_t\bullet D_s = \sum_w m(t,s,w) D_w,
\]
where the sum is over all rooted trees. For our aims it is better to
use a different description of this pre-Lie algebra. We closely follow
the strategy of Hoffman \cite{Hoff:tree} in this respect. If
$m(t,s,w)\neq 0$, it is exactly the order of the orbit of the root of
the subtree $s$ under the action of 
the group $\Aut(w)$ of automorphisms of $w$. If $s$ and $t$ are
$n$-coloured trees and $v$ is a vertex in $t$, denote by
$t\circ_{(v,i)}s$ the tree obtained from $t$ and $s$ by connecting the
root of $s$ to the vertex $v$ by an edge of colour $i$. Let $n(t,s,w)$
be the number if vertices $v\in t$ such that $t\circ_{(v,i)} s=w$ for
some $i\in \mathbf{p}$. Then for any such vertex $v$, the order of the
orbit of $v$ in $t$ under the action of $\Aut(t)$ is exactly
$n(t,s,w)$. 

Define an other pre-Lie algebra structure $\bullet'$ on the same
vector space $P(C_n^*)$, by
\[
D_t\bullet'D_s = \sum_w n(t,s,w)\,D_w,
\]
and denote this pre-Lie algebra by $P(C_n^*)'$.
For a subtree $s\subset w$, denote by $\Aut^s(w)$ the
automorphisms of $w$ that pointwise fix $s$. Then, if $m(s,t,w)\neq 0$
we can write, following Hoffman \cite{Hoff:tree},
\[
\begin{split}
m(t,s,w) &= \frac{|\Aut(w)|}{|\Aut^s(w)|\cdot|\Aut(t)|}\\
n(t,s,w) &= \frac{|\Aut(s)|}{|\Aut^{\{v\}}(t)|},
\end{split}
\]
for a vertex $v$ such that $t\circ_{(v,i)}s = w$ for some $i$.
Since $|\Aut^s(w)|= |\Aut^{\{v\}}(t)|$ it follows that $D_t\longmapsto |\Aut(t)|D_t$
defines an isomorphism of pre-Lie algebras $P(C_n^*)\longrightarrow
P(C_n^*)'$ (in characteristic 0).

In the remainder of the proof, let $\mathbf{p}=\{1,2,\ldots,n\}$. 
We prove \textit{(ii)} by constructing an inclusion
$P(C_n^*)'\longrightarrow L_n$. Note that
\[
D_t\bullet'D_s = \sum_{v\in t}\sum_{i\in \mathbf{p}} D_{t\circ_{(v,i)} s}
\]
For an $n$-coloured tree $t$, denote by ${\uparrow}_i(t)$ the tree with
coloured vertices obtined by moving the colour of ech edge up to the
vertex directly above it and colouring the root by $i$. Note that 
\[
{\uparrow}_j(t\circ_{(v,i)}s) = {\uparrow}_i(t) \circ_v {\uparrow}_j(s). 
\]
Let $\mathbf{p} = \{1,\ldots, n\}$ and consider $S(T_n)$ with the
corresponding Hopf algebra structure as defined above. 
Define $\phi:P(S(T_n)^*)'\longrightarrow L_n$ from the pre-Lie algebra
of primitives to the free pre-Lie algebra on $n$ generators by
\[
\phi(D_t) = \sum_{j=1}^n {\uparrow}_j(t).
\]
Then $\phi$ is a linear enbedding. Moreover, $\phi$ preserves the
pre-Lie algebra structure since
\[\begin{split}
\phi(D_t\bullet' D_s) &= \sum_{(v,i)} \phi(D_{t\circ_{(v,i)}s}) \\
&= \sum_{(v,i)} \sum_j {\uparrow}( t\circ_{(v,i)} s)\\
&= \sum_v \sum_{i,j} {\uparrow}_i(t) \circ_v {\uparrow}_j(s)\\
&= \phi(D_t) \bullet \phi(D_s).
\end{split}
\]
\end{Pf}

\begin{Ex}
In the case of the Connes-Kreimer Hopf algebra (Example
\ref{Ex:Connes-Kreimer}),
Theorem \ref{Tm:preLie} states that the dual Hopf algebra is the
universal enveloping algebra of the free pre-Lie algebra on one
generator as first proved by Chapton-Livernet \cite{ChapLiv:Prelie}. 
\end{Ex} 

\section{Formal coalgebra deformations}

\label{deform}
In this section we consider how the different Hopf algebra structures
of Theorem \ref{Cr:Hopfspelledout} are related from the point of view
of deformation theory. For simplicity we restrict our attention to the case $n=1$.  

Let $A$ be a Hopf algebra. Recall (e.g. Gerstenhaber-Shack \cite{GerSch:Bidef})
the bicomplex $C^{p\,q}(A)=\mathrm{Hom}(A^{\otimes
  p},A^{\otimes q})$ for $p,q\geq 1$. For $q$ fixed it is the
Hochschild complex of the algebra $A$ with coefficients in $A^{\otimes
  q}$, and for $p$ fixed it is the Hochschild complex of the coalgebra
  $A$ with coefficients in $A^{\otimes p}$. Classes in
  $H^3(\mathrm{Tot}(C^{**}(A)))$  are in 1-1 correspondence with
  Hopf algebra deformations of $A$ modulo $\hbar^2$.

A formal coalgebra deformation
$\Delta_\hbar$ of a Hopf algebra $A$ is a $k[[\hbar]]$-linear map 
$\Delta_\hbar:A[[\hbar]]\longrightarrow A[[\hbar]]\otimes_{k[[\hbar]]}A[[\hbar]$ for
which $\Delta_\hbar$ makes $A[[\hbar]]$ with the same multiplication
and counit a Hopf algebra
over $k[[\hbar]]$, and such that evaluation at $\hbar=0$ gives the original
Hopf algebra structure on $A$. 

The vector spaces $\text{Der}(A,A^{\otimes q})$ of algebra derivations
form the kernel of the horizontal differential at the edge of the
complex $C^{*q}(A)$, and thus a subcomplex of the coalgebra Hochschild
complex. Coalgebra deformations modulo $\hbar^2$ of the Hopf
$P$-algebra $A$ are in 1-1 correspondence with classes in
$H^2(\text{Der}(A,A^{\otimes *}))$. 

Let us now turn to the example of the Hopf algebra 
$S(T_1)$ with the coproduct on trees given by  
$\Delta(s) = s\otimes 1 + 1\otimes s$. This is the coproduct induced
by $\sigma_1=\sigma_2 = u \circ \hbar$.
If $q_1,q_2\in t\cdot k[[\hbar]]$ and if we write
$\lambda:=\lambda_1$, then the map
$\Delta_{q_1,q_2}:S(T_1)[[\hbar]]\longrightarrow S(T_1)[[\hbar]]
\otimes_{k[[\hbar]]} S(T_1)[[\hbar]]$ inductively defined by 
\[ 
\Delta_{q_1,q_2}\circ\lambda(x) = \sum_{(x)}\lambda(x_1)\otimes q_2^{|x_2|}x_2 +
q_1^{|x_1|}x_1\otimes \lambda(x_2)
\]
defines a coalgebra deformation of $S(T_1)$.

The Hopf algebra $S(T_1)$ is graded. Let us write
$\text{Der}_0(S(T_1),S(T_1)^{\otimes *})$ for the subcomplex  of
$\text{Der}(S(T_1),S(T_1)^{\otimes *})$ consisting of those 
derivations that preserve the degree. Classes in
$H^2(\text{Der}_0(S(T_1),S(T_1)^{\otimes *}))$ correspond to graded
coalgebra deformations. The result below studies the deformations
$\Delta_{q_1,q_2}$ as graded coalgebra deformations. 
%

\begin{Pp}\label{Lm:dDer1}
Consider $S(T_1)$ with the coproduct induced by $\sigma_1=\sigma_2 =
u\circ\hbar$.
\begin{enumerate}
\item
The boundaries in $\text{Der}_0(S(T_1),S(T_1)^{\otimes 2})$ are the
derivations $\phi$ that can be written
as  
\[
\phi(s) = \sum_w c_{s,w}\bar{\Delta}(w),
\]
for all trees $s$, and constants $c_{s,w}\in k$, and where the sum
ranges over all forests $w$ such that $|s| = |w|$. As usual,
$\bar{\Delta} = \Delta - (\id\otimes 1 + 1\otimes \id)$.
\item
Let $q_1 \equiv c_1 \hbar,\ q_2 \equiv c_2 \hbar$, and $q_1' \equiv d_1 \hbar,
q_2' \equiv d_2 \hbar$ modulo $\hbar^2$. The
two graded coalgebra deformations $\Delta_{q_1,q_2}$ and
$\Delta_{q'_1,q'_2}$ are equivalent modulo $\hbar^2$ iff $c_1-c_2 =
d_1-d_2$. 
\end{enumerate}
\end{Pp}
\begin{Pf}
Let $\psi\in \text{Der}_0(S(T_1),S(T_1))$. Then $\psi(1) = 0$, and $\psi$ is
determined by its values on trees. Write $\psi$ in matrix form as
$\psi(s) = \sum_w c_{s,w} w$, where the sum ranges over forests, and
$c_{s,w}\in k$ are constants.  Since we assume $\psi$ is graded,
$c_{s,w} = 0$ if $|w|\neq |s|$.
Compute
\[
\begin{split}
d\psi(s)  &= (\psi(s)\otimes 1 + 1\otimes \psi(s)) - \Delta(\psi(s)) \\
&= \sum_w c_{s,w}  (w\otimes 1 + 1\otimes w) - c_{s,w} \Delta(w) \\
&= - \sum_wc_{s,w}\bar{\Delta} w.
\end{split}
\]
%
For $\psi:S(T_1)\longrightarrow S(T_1)$ as above define the endomorphism
$\Psi$ of $S(T_1)[[\hbar]]$ by $\Psi(x) = x+ \psi(x)t$ for $x\in S(T_1)[[\hbar]]$. Two
coalgebra deformations $\Delta_{q_1q_2}$ and $\Delta_{q_1'q_2'}$
are equivalent modulo $\hbar^2$ iff we can find a derivation $\psi$ such
that the corresponding map $\Psi$ satisfies for all $s$
\[
\begin{split}
\Delta_{q_1'q_2'}\circ\Psi(s) &\equiv
\Psi\otimes\Psi\circ\Delta_{q_1,q_2}(s), \qquad\text{or equivalently}\\ 
\Delta_{q_1'q_2'}(s+\psi(s)\hbar) &\equiv \Delta_{q_1,q_2}(s)+ (\psi\otimes
1)(\Delta_{q_1,q_2}(s))\hbar + (1\otimes \psi)(\Delta_{q_1,q_2}(s))\hbar
\end{split}
\]
modulo $\hbar^2$.
To compute this we only have to consider terms in
$\Delta_{q_1'q_2'}(s)$ and 
$\Delta_{q_1,q_2}(s)$ corresponding to $u\subset s$, such that
\[
\sum_{v\in u}p(v,u,s) + \sum_{v\in u^c}p(v,u^c,s)\leq 1,
\]  
which means that either $u=s$, or  $u^c=s$, or $s= u\circ_r\lambda(1)$,
or $s= u^c\circ_r\lambda(1)$. We only need to consider terms
in $\Delta_{q_1'q_2'}(\psi(x))$ corresponding to $u\subset s$, such that
\[
\sum_{v\in u}p(v,u,s) + \sum_{v\in u^c}p(v,u^c,s) = 0,
\]
which means
that either $u=s$, or $u^c = s$.  From this it follows that $\Psi$
defines an equivalence mod $\hbar^2$ of the two deformations iff the
degree 1 terms in $t$ match, which is to say
\[
\begin{split}
\sum_{\{u\subset s|s= \lambda(1)\circ_ru\}} d_1 u\otimes \lambda(1) + d_2
\lambda(1)\otimes u +\sum_w c_{s,w} \Delta(w)
= 
\\ 
\sum_{\{u\subset s | s= \lambda(1)\circ_ru\}}
c_1 u\otimes \lambda(1) + c_2
\lambda(1)\otimes u +c_{s,w} w\otimes 1 + c_{s,w} 1\otimes w.
\end{split}
\]
Of course $c_{s,w} w\otimes 1 + c_{s,w} 1\otimes w$ is the primitive
part of $\sum_w c_{s,w} \Delta(w)$.
For the equality we thus need  
$c_{s,w} = (c_1-d_1)m = (c_2-d_2)m_{su}$ if $w = \lambda(1)u$ for 
an $u$ such that $s = \lambda(1)\circ_r u$ and $m_{su}$ is cardinality 
of the orbit of the vertex $s-u$ under the automorphism group of $s$.
Choose the other $c_{s,w}$ equal to 0.
\end{Pf}

\section{The general approach in the associative case}
\label{asscase}
The construction of bialgebras can be performed in much greater
generality. Starting from an
arbitrary Hopf operad $P$ one  
can construct the operad $P[\lambda_n]$ which has as algebras
$P$-algebras together with an $n$-ary operation. Under conditions 
similar to those in Theorem  \ref{Cr:Hopfspelledout} one can find a
Hopf-$P$ algebra structure on the initial
$P[\lambda_n]$-algebra (see \cite{Moer:Hoa},
\cite{Pep:thesis}). The explicit calculations for $P=\mathrm{Com}_*$,
the operad of unital commutative algebras is presented in the previous
sections. These can also be done for the operad $\mathrm{Ass}_*$ of
unital associative algebras. In this section we briefly state some of the
results. 
\begin{Df}
Let the category $\mathcal{A}_n$ of unital associative algebras with an
$n$-ary map be the category which has as as objects pairs
$(A,\alpha)$ of an associative unital algebra $A$ and a linear map
$\alpha:A^{\otimes n} \longrightarrow A$, and which has as maps
$\mathcal{A}_n((A,\alpha),(B,\beta))$ the algebra homomorphisms
$f:A\longrightarrow B$ such that $\beta\circ f^{\otimes n} = f\circ
\alpha$. Let $(A_n,\lambda_n)$ be the initial object in the category
$\mathcal{A}_n$. 
Then $A_n$ can be described as the free associative
algebra on trees with edges coloured 
by $\{1,\ldots,n\}$ and at each vertex the incoming edges of each colour
endowed with a separate linear ordering (\note{planar $n$-trees}, for
short).
\end{Df}
\begin{Tm}
The tensor algebra $A_n$ on planar rooted trees with $n$-coloured edges
has a natural family of graded connected Hopf algebra structures,
indexed by sequences
$(q_{11},\ldots,q_{1n},q_{21},\ldots,q_{2n})\in k^{2n}$.
The inductive description of the coproduct is given by the formula of
Theorem \ref{Cr:Hopfspelledout}.
\end{Tm}
The $A_n$ together again form a simplicial algebra (cf. Remark
\ref{Rm:simpl}). 
\begin{Cr} \label{Rm:Assdeform}
Let $A_n$ be the free associative algebra on the planar rooted trees.
\begin{enumerate}
\item
The coproducts on $A_n$ of Corollary \ref{Cr:Hopfspelledout} are given
by the closed formula 
\[
\Delta(t) = \sum_{s\subset t}
\prod_j \left(\prod_{v\in s} q_{1\,j}^{p_j(v,s,t)}\cdot
\prod_{v\in s^c}q_{2\,j}^{p_j(v,s^c,t)}
s\otimes s^c\right),
\]
where the product of trees is associative, non-commutative. The order of
multiplication is given by the linear order on the roots of the trees 
defined by the linear on the incoming edges at each vertex and the partial 
order on vertices.
\item
The vector space of primitive elements of the dual is
spanned by elements dual to planar $n$-trees. The Lie bracket is the
commutator of the (non-associative) product $\bullet$ given by
\[
D_s\bullet D_t = \sum_w\sum_{s\subset w,\ s^c=t} 
\prod_j \left(\prod_{v\in s} q_{1\,j}^{p_j(v,s,w)}\cdot
\prod_{v\in t}q_{2\,j}^{p_j(v,t,w)}
D_w\right),
\]
where $w$, $s$ and $t$ are trees with a linear ordering on the incoming 
edges of the same colour at each vertex and the inclusions of $s$ and
$t$ in $w$ have to respect these orderings.
\end{enumerate}
\end{Cr}
\begin{Pf}
The only change is that we have to remember the ordering 
of up going edges at each vertex. Than one can copy the proof of
the commutative case verbatim.
\end{Pf}
\begin{Rm}
Independently, Foissy \cite{Fois:PlantreeII} has found the formula 
for the Lie bracket in  the case where  $n=1,\ q_{11}=1$ and $q_{12}=0$ (and 
$\text{Ass}_*$ is the underlying operad). He uses this formula to
give an explicit isomorphism between the Hopf algebras $A_1$ and
$A_1^*$ with this coproduct.
\end{Rm}

\bibliographystyle{plain}
\bibliography{hopf}

\begin{thebibliography}{10}

\bibitem{BrouFrab:Trees}
C.~Brouder and A.~Frabetti.
\newblock Renormalization of {QED} with {P}lanar {B}inary {T}rees.
\newblock {\em Euro. Phys. J. C}, 19(4):715--741, 2001.

\bibitem{ChapLiv:Prelie}
F.~Chapoton and M.~Livernet.
\newblock Pre-{L}ie algebras and the rooted trees operad.
\newblock {\em Internat. Math. Res. Notices}, 8:395--408, 2001.
\newblock {Also at \texttt{arXiv:math.QA/0002069}}.

\bibitem{ConKr:Hoa}
A.~Connes and D.~Kreimer.
\newblock Hopf {A}lgebras, {R}enormalisation and {N}oncommutative {G}eometry.
\newblock {\em Comm. Math. Phys.}, 199, 1998.

\bibitem{Fois:PlantreeII}
L.~Foissy.
\newblock Les alg\`ebres de {H}opf des arbres enracin\'es d\'ecor\'es. {II}.
\newblock {\em Bull. Sci. Math.}, 126(4):249--288, 2002.
\newblock Also at \texttt{arXiv:math.QA/0105212}.

\bibitem{GerSch:Bidef}
M.~Gerstenhaber and S.D. Schack.
\newblock Bialgebra cohomology, deformations, and quantum groups.
\newblock {\em Proc. Nat. Acad. Sci. U.S.A.}, 87(1):478--481, 1990.

\bibitem{Hoff:tree}
M.E. Hoffman.
\newblock Combinatorics of rooted trees and {H}opf algebras.
\newblock {\em Trans. Amer. Math. Soc.}, 355(9):3795--3811 (electronic), 2003.

\bibitem{Krei:Hoa}
D.~Kreimer.
\newblock On the {H}opf algebra structure of perturbative quantum field
  theories.
\newblock {\em Adv. Theor. Math. Phys.}, 2(2):303--334, 1998.

\bibitem{LodRon:Trees}
J.-L. Loday and Maria~O. Ronco.
\newblock Hopf {A}lgebra of the {P}lanar {B}inary {T}rees.
\newblock {\em Adv. Math.}, 139(2):293--309, 1998.

\bibitem{MilnorMoore:Hoa}
J.~Milnor and J.~Moore.
\newblock On the {S}tructure of {H}opf {A}lgebras.
\newblock {\em Ann. Math}, 81(2):211--264, 1965.

\bibitem{Moer:Hoa}
I.~Moerdijk.
\newblock On the {C}onnes-{K}reimer {C}onstruction of {H}opf {A}lgebras.
\newblock {\em Cont. Math.}, 271:311--321, 2001.
\newblock Also at {\texttt{arXiv:math-ph/9907010}}.

\bibitem{Pep:thesis}
{P.P.I. van der Laan}.
\newblock Operads - {H}opf algebras and coloured {K}oszul duality.
\newblock Ph.D. thesis Utrecht University, 2004.

\end{thebibliography}
~\\
Pepijn van der Laan (\texttt{vdlaan@math.uu.nl})\\ 
\\~\\
Ieke Moerdijk (\texttt{moerdijk@math.uu.nl})\\ 
Mathematisch Instituut, Universiteit Utrecht
\end{document}